\documentclass[twoside]{amsart}
\usepackage{graphicx}
\usepackage[usenames]{color}
\usepackage{amsthm, amsmath, amssymb, amsfonts, graphicx, epsfig, mathtools}
\usepackage{dsfont}
\usepackage{lscape}
\usepackage{cancel}
\usepackage[normalem]{ulem}
\usepackage[T1]{fontenc}

\newtheorem{thm}{Theorem}[section]
\newtheorem{lem}{Lemma}[section]
\newtheorem{pro}{Proposition}[section]
\newtheorem{cor}{Corollary}[section]
\newtheorem{axiom}{Definition}[section]
\newtheorem{rem}{Remark}[section]

\newtheorem{ex}{Example}[section]

\newcommand{\be}{\begin{equation}}
\newcommand{\ee}{\end{equation}}
\newcommand{\bde}{\begin{displaymath}}
\newcommand{\ede}{\end{displaymath}}
\newcommand{\beq}{\begin{eqnarray*}}
\newcommand{\eeq}{\end{eqnarray*}}
\newcommand{\beqa}{\begin{eqnarray}}
\newcommand{\eeqa}{\end{eqnarray}}
\newcommand{\bel }{\left\{\begin{array}{ll}}
\newcommand{\eel}{\cr \end{array} \right.}
\newcommand{\bd}{\begin{axiom} \rm }
\newcommand{\ed}{\end{axiom} \rm }
 \newcommand{\bex}{\begin{ex} \rm }
\newcommand{\eex}{\end{ex}}
\newcommand{\bt}{\begin{thm}}
\newcommand{\et}{\end{thm}}
\newcommand{\bl}{\begin{lem}}
\newcommand{\el}{\end{lem}}
\newcommand{\bp}{\begin{pro}}
\newcommand{\ep}{\end{pro}}
\newcommand{\bcor}{\begin{cor}}
\newcommand{\ecor}{\end{cor}}
\newcommand{\lab }{\label }
\newcommand{\br}{\begin{rem}}
\newcommand{\er}{\end{rem}}

\def\proof{\noindent {\it Proof. $\, $}}
\def\finproof{\hfill $\Box$ \vskip 5 pt}

\def\I{\mathds{1}}

\def \wh{\widehat}

\def\r{\mathbb R}
\def\F{{\mathcal F}}

\def\ff{{\mathbb F}}

\def\P{\mathbb P}

\def\EP{{\mathbb E}_{\mathbb P}}

\newcommand{\set}[1]{{\left\{#1\right\}}}
\newcommand{\Real}{\mathbb R}

\def\wrt{\mbox{with respect to}}
 \def\KS{{K_1}}
\def\KO{{K_2}}

    \usepackage{fancyhdr}
    \usepackage{lastpage}
\title{{\Large \bf
 Intricacies of Dependence between Components of Multivariate
Markov Chains: Weak Markov Consistency and Weak Markov Copulae
} \vskip 80 pt}
\author{Tomasz R.  Bielecki}\thanks{Research of T.R. Bielecki was partially supported by NSF Grant DMS-0604789 and NSF Grant DMS-0908099.}
\address{
Department of Applied Mathematics \\
Illinois Institute of Technology \\
Chicago, IL 60616, USA \\
E-mal: {\tt bielecki@iit.edu }
}
\author{Jacek Jakubowski}\thanks{Research of J. Jakubowski and M. Niew\k{e}g\l owski was partially supported by Polish MNiSW grant N N201 547838.}
\address{
Institute of Mathematics, University
of Warsaw \\
 ul. Banacha 2, 02-097 Warszawa, Poland \\
 and \\
 Faculty of Mathematics and Information Science  \\
Warsaw University of
Technology\\
ul. Koszykowa 75, 00-662 Warszawa, Poland \\
E-mail: {\tt jakub@mimuw.edu.pl }
}
\author{Mariusz Niew\k{e}g\L owski}
\address{
Faculty of Mathematics and
Information Science, Warsaw University of Technology  \\
ul. Koszykowa 75, 00-662 Warszawa, Poland \\
E-mail:  {\tt
M.Nieweglowski@mini.pw.edu.pl}
}

\begin{document}
\begin{abstract}
In this paper we examine the problem of existence and construction of multivariate Markov chains such that their components are Markov chains with given laws. Specifically, we provide sufficient and necessary conditions, in terms of semimartingale characteristics, for a component of a multivariate Markov chain to be a Markov chain in its own filtration - a property called weak Markov consistency. Accordingly, we introduce and discuss the concept of weak Markov copulae. Finally, we examine relationship between the concepts of weak Markov consistency and weak Markov copulae, and the corresponding strong versions of these concepts.

\end{abstract}
\maketitle
\begin{center}
First version circulated since 13 May 2011
\end{center}

\tableofcontents

\section*{Introduction}

Modeling of dependence between stochastic processes is a very
important issue arising from many different  applications, among
others in financial mathematics.   By modeling dependence we mean
construction of  a multivariate stochastic process with prescribed
marginal laws. In this paper we focus on Markov chains,  and deal
with the problem of constructing a multivariate Markov chain such
that its components are given Markov chains in their own filtrations. It is
well known that components of multivariate Markov process  are in
general not Markovian (in any filtration), so the problem that we study here is by no means a trivial one.
We give sufficient and necessary conditions, in terms of the
semimartingale characteristics, for a component of a multivariate
Markov chain to be a Markov chain in its own filtration.

%
Our paper continues the study of Markovian consistency and Markov copulae  for
multivariate Markov processes, initiated in \cite{BVV}, \cite{BJVV},
\cite{BJN2} and \cite{BJN}.

 Here, we introduce and study the concept of {\it weak Markovian
consistency}, and we relate it to the concept of {\it strong
Markovian consistency} that was explored in the aforementioned
papers under the name of {\it Markovian consistency}.
 We also continue the study of dependence between Markov processes.
Thus, we continue the study of Markov copulae, the concept
originally introduced in \cite{BVV}. Specifically, we introduce and examine Markov copulae
with regard to weak Markovian consistency. It turns out that certain
unwanted features of Markov copulae, inherent to the framework of
strong Markovian consistency, are no longer present in the framework
of weak Markovian consistency. This is particularly pleasing in view
of applications of Markov copulae in credit risk management or in reliability management; in fact, this aspect of weak Markov
copulae makes them exceptionally important tool in modeling  dynamic dependence. We provide more insight into this important issue in  Remark \ref{fin}.

 As already said, we confine our discussion, for the most part, to the case of finite Markov chains. One might object the choice of finite Markov chains as the of object of interest in this paper, as  one might think that this choice is very restrictive. In \cite{BJN} we studied strong Markovian dependence in the context of (nice) Feller processes. What we learned while working on paper \cite{BJN} and while working on the present paper,  is that from the point of view of intricacies of dependence between components of a multivariate Markov process, the finite state space set-up is actually not restrictive at all! The dependence here is equally intricate as dependence in the case of general Feller process, which is much harder to present, due to various technicalities that obscure the dependence picture. That is why, with the benefit for the reader, we are presenting here a study of the intricate dependence between components of finite Markov chains, which does not require any use of sophisticated technical machinery, but at the same allows for pointing to the essence of the of intricacies of dependence between components of a multivariate Markov process.

It needs to be noted that problems that we study in the present
paper are also connected with lumpability problem for continuous
time Markov chains  (see Ball and Yeo \cite{BY93} and discussion
there, Burke and Rosenblatt \cite{BR58}). In \cite{BY93} necessary
and sufficient conditions are provided for intensity matrix so that
the marginal component process of a Markov chain is a time
homogenous continuous time Markov chain in its natural filtration.
If we omit the assumption of time homogeneity and weaken assumption
on intensity matrix, then there exist Markov process with marginals
being also Markov in their own filtration which does not satisfy
conditions from \cite{BY93} (see Example 3.2.). Moreover assumptions
imposed in these papers on intensity matrix exclude Markov chains
with absorbing states, a case that can be treated using our
methodology.

 In case of a bivariate Markov
chain, one can also note some similarity between our work and the
studies of Markovian coupling (see e.g.  Chen \cite[Sect.
I.5.2]{Chen2004}). It needs to be stressed though that the
concepts of weak Markovian consistency and weak Markovian copulae
are much more than (standard) coupling of  Markov chains; and this
not just because these concepts apply to multivariate case and not
only to the bivariate case.  Markovian couplings deal with
marginal properties of transition probabilities (cf. equations
(5.7) in \cite{Chen2004}) and with ``marginal'' properties of
generators (cf. equations (5.8) in \cite{Chen2004}). Specifically,
the properties looked at within the Markovian coupling universe,
that are somewhat relevant to our present work,  amount to
\begin{enumerate}
\item the property that marginals of a bivariate transition probability are equal to given univariate transition probabilities (cf. equations (5.7) in \cite{Chen2004}), and
\item
the property that ``univariate projections'' of a bivariate
Markovian generator are equal to given univariate Markovian
generators (cf. equations (5.8) in \cite{Chen2004}).
\end{enumerate}
This however is much less than dealing with the marginal laws of a process in the sense of
\begin{itemize}
\item
asking questions regarding Markovian consistency: that is, asking
questions regarding necessary and sufficient conditions that need
to be satisfied by the generator of a multivariate Markov chain,
such that the chain's components are Markovian either in their own
filtrations (the property of weak Markovian consistency), or are
Markovian either in the filtration of the entire multivariate
process (the property of strong Markovian consistency), or are not
Markovian at all.
\item
asking questions regarding  Markov copulae: that is, asking
questions regarding construction of a generator of a multivariate
Markov chain, such that the chain's components are Markovian either
in their own filtrations and their laws coincide with the laws of
given univariate Markov chains (i.e. construction of a weak
Markovian  copula), or the chain's
components are Markovian in the filtration of the entire
multivariate process and their laws coincide with the laws of
given univariate Markov chains (i.e. construction of a strong
Markovian  copula).
\end{itemize}
In addition, studies of Markovian coupling do not touch the issues of intricate nature of
dependence between components of a multivariate Markov chain, that we study and demonstrate in this paper.

The paper is organized as follows. In Section 1 we give a sufficient
and necessary condition for a multivariate Markov chain to be weakly
consistent. Note that a sufficient condition for weak Markovian
consistency can be deduced from the result of Rogers and Pitman
\cite{R-P} in which   sufficient conditions for  a function of a
Markov process to be  a Markov process are given. Our condition for a
weak Markovian consistency is not only more explicit, but also
necessary. We also study the question when weak Markovian
consistency implies strong Markovian consistency. It turns out that
this is equivalent to $\P$-immersion between $\ff^{X^i}$ and
$\ff^X$, given that weak Markovian consistency holds.
 In Section 2 we  study weak Markov copulae.
In Section 3  we present  three simple, but non-trivial examples,
that illustrate \textbf{intricacies of dependence} between components of a
multivariate Markov chain. Specifically, in Examples 3.1--3.3 we show that
\begin{enumerate}
\item there exist Markov processes that are strongly Markovian consistent,
\item there exist Markov processes that are weakly Markovian consistent, but are not strongly Markovian consistent; in addition, in this case, one would expect that even if a multivariate Markov process is time-homogeneous, its components are time-inhomogeneous Markov processes; Example 3.2 illustrates this,
\item there exist Markov processes that are neither strongly Markovian consistent nor weakly Markovian consistent.
\end{enumerate}

\section{Markovian Consistency}

As already said, we shall focus in this paper on the case of finite Markov chains. Nevertheless, we shall formulate  the concept of weak Markovian consistency in more generality. Towards this end we consider $X=(X^n,\, n=1,\dots,N)$, a multivariate Markov process, defined on an underlying probability space $(\Omega,
{\mathcal F}, \mathbb{P})$, taking values in $\Real^N$.\footnote{The study presented in this paper carries over to the case of multivariate Markov process taking values in a product of arbitrary (metric) spaces.} We denote by $\ff^X$ the filtration of $X$, and by $\ff^{X^n}$ the filtration of the coordinate $X^n$ of $X.$
It is well known that, in general, the coordinates of $X$ are not Markov with respect to their own filtrations.

\bd (i) Let us fix $n$. We say that the process $X$ satisfies the
\textit{weak Markovian consistency condition with respect to the
component $X^n$} if  for every $B \in \mathcal{B} (\mathbb{R}) $
and all $t,s \geq 0 $,
\begin{equation}\label{eq:equal-in-trans-prob-1}
\mathbb{P}\left( X^n_{t + s } \in B | \mathcal{F}^{X^n}_t \right) = \mathbb{P} \left( X^n_{t + s } \in B | X^n_t   \right),
\end{equation}
so that the component $X^n$ of $X$ is a Markov process in its own filtration. \\
(ii) If $X$ satisfies the weak Markovian consistency condition
with respect to $X^n$ for each $n \in \set{ 1,\ldots,N }$, then we say that
$X$ satisfies the weak Markovian consistency condition.
\ed

Previously, in  \cite{BVV}, \cite{BJVV}, \cite{BJN2} and \cite{BJN}, a stronger concept was studied.

\bd (i) Let us fix $n$. We say that the process $X$ satisfies the
\textit{strong Markovian consistency condition with respect to the
component $X^n$} if  for every $B \in \mathcal{B} (\mathbb{R}) $
and all $t,s \geq 0$,
\begin{equation}\label{eq:markov-cons}
\mathbb{P}\left( X^n_{t + s } \in B | \mathcal{F}^X_t \right) =   \mathbb{P} \left( X^n_{t + s } \in B | X^n_t    \right)
\end{equation}
or equivalently
\begin{equation}\label{eq:equal-in-trans-prob}
\mathbb{P}\left( X^n_{t + s } \in B | X_t  \right) =  \mathbb{P} \left( X^n_{t + s } \in B | X^n_t    \right),
\end{equation}
so that $X^n$ is a Markov process in the filtration of $X$.\\
(ii) If $X$ satisfies the strong Markovian consistency condition with respect to $X^n$ for each $n \in \set{ 1,\ldots,N }$, then we say that $X$ satisfies the strong Markovian consistency condition.
\ed

Obviously, strong Markovian consistency implies weak Markovian consistency, but not vice versa as will be seen in one of the examples in Section \ref{examples}. As a matter of fact, it may happen that all components of $X$ are Markovian in their filtrations, but $X$ is not Markovian in its filtration (see e.g. Bielecki et al. \cite[Example 2.4.2]{BJN2}).


From now on we assume that $X=(X^1,\ldots,X^N)$ is a Markov chain
with values in a finite product space, say ${\mathcal X}={\sf
X}_{n=1}^N {\mathcal X}^n,$ where ${\mathcal
X}^n=\{x^n_1,\ldots,x^n_{m_n}\}\subseteq \r.$ However, to somewhat
simplify the notation, in most of the paper we shall consider
bivariate processes $X$ only, that is, we put $N=2$, and we take
$\Lambda(t) =[\lambda^x_y(t)]_{x,y\in {\mathcal X}}$ as a generic
symbol for the $\P$-infinitesimal generator of $X$. Thus, $\Lambda
(t)$ is an $m\times m$ matrix, where $m=m_1\cdot m_2.$ We stress
that restriction to bivariate case is for a notational convenience
only. Our results naturally extend to the multivariate case.
%

\subsection{Semimartingale characterization of a finite Markov chain}
Let us consider a 
c\`{a}dl\`{a}g process $V$ defined on   $(\Omega,
{\mathcal F}, \mathbb{P})$, taking values in a finite set ${\mathcal V} \subset \r^N$.


For any two distinct states $v,w\in {\mathcal V},\, $ we
define an $\ff^V$-optional random measure $N_{vw}$ on
$[0,\infty)$ by
\be\lab{N} N_{vw}((0,t])=\sum_{0<s\leq
t}\I_{\{V_{s-}=v,V_{s}=w\}}.\ee \noindent We shall simply write
$N_{vw}(t)$ in place of $N_{vw}((0,t]).$ Manifestly, $N_{vw}(t)$
represents the number of jumps from state $v$ to state $w$ that the
process $V$ executes over the time interval $(0,t].$ Let us denote
by $\nu_{vw}$ the dual predictable projection (the \emph{compensator}) with respect to  $\ff^V$ of the
random measure $N_{vw}$.


Next, let us define a deterministic matrix valued function $\Lambda$ on $[0,\infty)$
by \be \Lambda (t)=[\lambda^v_w(t)]_{v,w\in {\mathcal V}}, \ee where
$\lambda^v_w$'s are real valued, locally integrable functions on
$[0,\infty)$ such that for $t\in [0,\infty )$ and $v,v\in {\mathcal
V},\, v\ne w$, we have
$$\lambda^v_w(t)\geq 0$$
\noindent and
$$\lambda^v_v(t)=-\sum_{w\ne v}\, \lambda^v_w(t).$$ \noindent

The following result, gives necessary and sufficient condition for c\`{a}dl\`{a}g process $V$ with values in $\mathcal{V}$ to be a Markov chain.

\bp\lab{prop:comp} A process $V$
is a Markov chain ($\wrt$ $\ff^V$) with infinitesimal generator
$\Lambda (t)$ iff the compensators
 with respect to  $\ff^V$ of the counting measures $N_{vw}(dt)$, $v,w\in
\mathcal{V}$, are of the form
 \be\lab{nu}
\nu_{vw}((0,t])=\int_0^t \I_{\{V_{s}=v\}}\lambda^v_w(s)ds.
\ee
\ep

\proof It has been shown in Lemma 5.1 in  \cite{BJVV} that a process $V$
is a Markov chain ($\wrt$ $\ff^V$) with infinitesimal generator
$\Lambda (t)$ iff the compensators
 with respect to  $\ff^V$ of the counting measures $N_{vw}(dt)$, $v,w\in
\mathcal{V}$, are of the form
 \be\lab{nu1}
\nu_{vw}((0,t])=\int_0^t \I_{\{V_{s-}=v\}}\lambda^v_w(s)ds.
\ee
Now, analysis of the proof of Lemma 5.1 in  \cite{BJVV} indicates that  the left hand limits $V_{t-}$ used in Lemma 5.1 in  \cite{BJVV} can, in fact, be replaced with $V_t$, which proves the present result.
\finproof

\begin{rem}
A finite Markov chain $V$ with a locally integrable generator $\Lambda(t)$ is a semimartingale (see, e.g., Elliott et al. \cite[Chapter 7.2]{EAM1994}). The jump measure of $V$, say $\mu^V$, can be expressed in terms of summation of the jump measures $N_{vw}$. Thus, in view of Proposition \ref{prop:comp} the infinitesimal characteristic of $V$ (with respect to an appropriate truncation function), which is the compensator of $\mu^V$ (denoted by $\nu^V$\!) is given in terms of summation of the compensators $\nu_{vw}$.
Indeed, one can easily check that if we define a truncation function $h$ by
\[
h(x) := x \I_\set{|x| \leq d },
\quad
\text{  where   }
\quad
     d := \frac{1}{2}\min \set{ |v -w| : v \neq w , v \in \mathcal{V}, w \in \mathcal{V} },
\]
then  $(0,0,\nu^V)$ is the local characteristic of $V$, where
\[
    \nu^V(dx,dt) = \sum_{v,w \in \mathcal{V}: v \neq w} \delta_{  w-v}(dx) \nu_{vw}(dt),
\]
and $\delta$ denotes the Dirac measure.
\end{rem}

\subsection{Necessary and sufficient conditions for weak Markovian consistency in terms of semimartingale characteristics}

Let us recall that we consider  bivariate processes. We take $n=1$ and we study the weak Markovian consistency of $X$ with respect to $X^1$. A completely analogous discussion can be carried out with respect to $X^2$.

For any two states $x^1,y^1\in {\mathcal X}^1$ such that $x^1\ne y^1,\, $ we
define the following $\ff^X$-optional random measure on
$[0,\infty)$:

\be\lab{N1} N^1_{x^1y^1}((0,t])=\sum_{0<s\leq
t}\I_{\{X^1_{s-}=x^1,X^1_{s}=y^1\}}.\ee

\noindent We shall write
$ N^1_{x^1y^1}(t)$ in place of $ N^1_{x^1y^1}((0,t]),$ and we shall denote
by $ \nu^1_{x^1y^1}$ the dual predictable projection (the compensator) with respect to  $\ff^X$ of the
random measure $ N^1_{x^1y^1}$.

Next, for any two states $x=(x^1,x^2), y=(y^1,y^2)\in {\mathcal X}$ such that $x\ne y,\, $ we
define  an $\ff^X$-optional random measure on
$[0,\infty)$ by

\begin{align}
\lab{N12} N_{xy}((0,t])=\sum_{0<s\leq
t}\I_{
\{
(X^1_{s-}=x^1,X^2_{s-}=x^2),(X^1_{s}=y^1,X^2_{s}=y^2)
\}
}.
\end{align}

\noindent We shall write
$ N_{xy}(t)$ in place of $ N_{xy}((0,t]),$ and we shall denote
by $ \nu_{xy}$ the compensator  of  $ N_{xy}$ with respect to  $\ff^X$.

It is easy to see that
\begin{equation}\label{eq:N1}
 N^1_{x^1y^1}(t)=\sum_{x^2,y^2\in {\mathcal X}^2}\, N_{(x^1,x^2),(y^1,y^2)}(t),
\end{equation}
and consequently (due to uniqueness of compensators)
\begin{equation}\label{eq:Nu1}
 \nu^1_{x^1y^1}((0, t])=\sum_{x^2,y^2\in {\mathcal X}^2}\, \nu_{(x^1,x^2),(y^1,y^2)}((0, t]).
\end{equation}
In view of Proposition \ref{prop:comp}, we see that for any two distinct states $x=(x^1,x^2), y=(y^1,y^2)\in {\mathcal X},\, $
\be\lab{nu12}
\nu_{(x^1,x^2),(y^1,y^2)}(dt)=\I_{\{(X^1_{t},X^2_t)=(x^1,x^2)\}}\lambda^{x^1x^2}_{y^1y^2}(t)dt. \ee
Let us denote by $\wh \nu^1_{x^1y^1}$ the compensator of the measure $ N^1_{x^1y^1}$  with respect to  $\ff^{X^1}$.
\bl\label{lem:komp_N_1}
 Assume that $X$ is a Markov chain with respect to its own filtration.
 The $\ff^{X^1}$\!\!-compensator  of  $N^1_{x^1,y^1}$ has the form
 \beqa\label{eq:komp_N_1}
    \wh \nu^1_{x^1y^1}(dt)
    =\I_{\{X^1_{t}=x^1\}}\sum_{x^2,y^2\in {\mathcal X}^2}\lambda^{x^1x^2}_{y^1y^2}(t)\EP (\I_{\{X^2_t=x^2\}}|{\mathcal F}^{X^1}_t)dt.
 \eeqa
\el
\begin{proof}
It follows from Lemma 4.3 in \cite{BJVV} that
\beqa
\label{eq:komp_N_1a}
\wh \nu^1_{x^1y^1}(dt) &=& \sum_{x^2,y^2\in {\mathcal X}^2}\EP (\I_{\{(X^1_{t},X^2_t)=(x^1,x^2)\}}\lambda^{x^1x^2}_{y^1y^2}(t)|{\mathcal F}^{X^1}_{t-})dt \\
&=& \nonumber
\sum_{x^2,y^2\in {\mathcal X}^2}\lambda^{x^1x^2}_{y^1y^2}(t)\EP (\I_{\{X^1_{t}=x^1\}}\I_{\{X^2_t=x^2\}}|{\mathcal F}^{X^1}_{t-})dt.
\eeqa
The process $X$  is quasi-left continuous, since it is a Markov chain. Hence, $X^1$ is also quasi-left continuous, so its natural filtration $\mathbb{F}^{X^1}$ is  quasi-left continuous and hence $\mathcal{F}^{X^1}_t = \mathcal{F}^{X^1}_{t-}$ (see Rogers and Williams \cite[III.11]{RW2000}). Thus by \eqref{eq:komp_N_1a} we have \eqref{eq:komp_N_1}.
\finproof
\end{proof}
Using Lemma  \ref{lem:komp_N_1} and Proposition \ref{prop:comp} we obtain the following important result.

\bt\lab{wmc} The component $X^1$ of $X$ is a Markov chain with respect to its own filtration if and only if
\begin{align}\lab{vn}
\I_{\{X^1_{t}=x^1\!\}}\!\!\!\!\sum_{x^2,y^2\in {\mathcal X}^2} \!\!\!\!\lambda^{x^1x^2}_{y^1y^2}(t)
\EP \left(\I_{\{X^2_t=x^2\}}|{\mathcal F}^{X^1}_t\!\right)
\!=\! \I_{\{X^1_{t}=x^1\}}\lambda^1_{x^1y^1}(t) \quad {dt \otimes
d\mathbb{P} \text{-a.s.} }  \  \forall x^1\!,y^1\!\in\! {\mathcal
X}^1\!,\ x^1 \neq y^1
\end{align}
for some locally integrable functions $\lambda^1_{x^1y^1}.$ The
generator of $X^1$ is $\Lambda^1 (t)
=[\lambda^1_{x^1y^1}(t)]_{x^1,y^1\in {\mathcal X}^1}$   with
$\lambda^1_{x^1x^1}$ given by
\[
    \lambda^1_{x^1x^1}(t) = - \sum_{y^1 \in \mathcal{X}^1, y^1 \neq x^1 } \lambda^1_{x^1y^1}(t)
    \qquad
    \forall x^1 \in \mathcal{X}^1
    .
\]
\et

\begin{proof}
Assume that \eqref{vn} holds. Since $X$ is a Markov chain, for
each $x^1, y^1 \in \mathcal{X}^1$, the $\mathbb{F}^{X^1}$
compensator of $N^1_{x^1,y^1}$ has, by Lemma \ref{lem:komp_N_1}
and \eqref{vn}, the form
\[
    \widehat{\nu}^1_{x^1y^1}(dt) = \I_\set{ X^1_t = x^1} \lambda^1_{x^1y^1}(t) dt
\]
for some locally integrable deterministic function $\lambda^1_{x^1y^1}$. In particular, note that \eqref{vn} implies that  $\lambda^1_{x^1y^1}$ is non-negative for   $x^1\ne y^1.$ Then, by Proposition \ref{prop:comp}, $X^1$ is a Markov chain with generator $\Lambda^1 (t) =[\lambda^1_{x^1y^1}(t)]_{x^1,y^1\in {\mathcal X}^1}.$
Conversely, assume that $X^1$ is a Markov chain with respect to its natural filtration with generator $\Lambda^1 (t) =[\lambda^1_{x^1y^1}(t)]_{x^1,y^1\in {\mathcal X}^1}$. Then \eqref{vn}  follows from Lemma \ref{lem:komp_N_1} and Proposition \ref{prop:comp}.
\finproof

\end{proof}

\begin{rem}
Note that \eqref{vn} implies that
\begin{align}\lab{vn1}
\I_{\{X^1_{t}=x^1\}}\!\!\!\!\!\!\sum_{x^2,y^2\in {\mathcal
X}^2}\!\!\!\!\!\!\lambda^{x^1x^2}_{y^1y^2}(t)\EP
\left(\I_{\{X^2_t=x^2\}}|X^1_t=x^1 \right)=
\I_{\{X^1_{t}=x^1\}}\lambda^1_{x^1y^1}(t) \quad {dt \otimes
d\mathbb{P} \text{-a.s.} } \  \forall x^1\!,y^1\!\in\! { \mathcal
X}^1\!,\ x^1     \neq y^1.
\end{align}
Thus, condition \eqref{vn1} is necessary for the weak Markovian consistency of $X$  with respect to  $X^1$.
\end{rem}

\subsection{Necessary and sufficient conditions for strong Markovian consistency }
Since one of our goals is to relate the notions of weak and strong Markovian consistency, we shall discuss in this section necessary and sufficient conditions for strong Markovian consistency of our finite Markov chain. Towards this end let us first recall condition (M) from \cite{BJVV}:

\noindent {\sf Condition (M):} {\it The generator matrix function $\Lambda$ satisfies for every $t \geq 0$
\be\tag{M1}\sum_{y^2\in \mathcal{X}^2}\lambda^{x^1x^2}_{y^1y^2}(t)=\sum_{y^2\in \mathcal{X}^2}\lambda^{x^1\bar x^2}_{y^1y^2}(t),
\  \, \forall x^2,\bar x^2 \in \mathcal{X}^2,\ \forall x^1,y^1\in
\mathcal{X}^1,\, x^1\ne y^1,
\ee and
\be\tag{M2}\sum_{y^1\in \mathcal{X}^1}\lambda^{x^1x^2}_{y^1y^2}(t)=\sum_{y^1\in \mathcal{X}^1}\lambda^{\bar x^1 x^2}_{y^1y^2}(t),
\  \, \forall x^1,\bar x^1 \in \mathcal{X}^1,\ \forall x^2,y^2\in
\mathcal{X}^2,\, x^2\ne y^2.
\ee
}
Next, consider the functions $\lambda^1_{x^1y^1}$ given, for $t \geq 0,$ by
\be\lab{lhat} \lambda^1_{x^1y^1}(t)=\sum_{y^2\in \mathcal{X}^2}\lambda^{x^1x^2}_{y^1y^2}(t),
\quad  x^1,y^1 \in \mathcal{X}^1,\, x^1\ne y^1,\quad
\lambda^1_{x^1x^1}(t)=-\sum_{y^1\in \mathcal{X}^1,y^1\ne x^1}\, \lambda^1_{x^1y^1}(t),\quad
\forall x^1\in \mathcal{X}^1. \ee

Under condition (M1), the functions ${\lambda}^1_{x^1y^1}$ are well defined and locally integrable, and it is straightforward to verify that they satisfy \eqref{vn}, so that weak Markovian consistency holds for $X$ with respect to  $X^1$.

Result analogous to Theorem \ref{wmc}, but with respect to component $X^2$, reads:
\begin{itemize}\item
\textit{The process $X^2$ is a Markov chain with respect to its own filtration if and only if
\begin{align}\lab{vn2}
\I_{\{X^2_{t}=x^2\}}\!\!\!\!\!\sum_{x^1,y^1\in {\mathcal
X}^1}\!\!\!\!\!\lambda^{x^1x^2}_{y^1y^2}(t)\EP
\left(\I_{\{X^1_t=x^1\}}|{\mathcal F}^{X^2}_t \!\right)=
\I_{\{X^2_{t}=x^2\}}\lambda^2_{x^2y^2}(t) \quad {dt \otimes
d\mathbb{P} \text{-a.s.} } \ \forall x^2\!,y^2\!\in\! {\mathcal
X}^1\!, x^2 \neq y^2
\end{align}
for some locally integrable functions $\lambda^2_{x^2y^2}.$ Then the generator of $X^2$ is $\Lambda^2 (t) =[\lambda^2_{x^2y^2}(t)]_{x^2,y^2\in {\mathcal X}^2}$
 with $\lambda^2_{x^2x^2}$ given by
\[
\lambda^2_{x^2x^2}(t)=-\sum_{y^2\in \mathcal{X}^2,y^2\ne x^2}\, \lambda^2_{x^2y^2}(t),\quad
\forall x^2\in \mathcal{X}^2.
\]}
\end{itemize}
Now, if we define
\be\lab{lhat2} \lambda^2_{x^2y^2}(t)=\sum_{y^1\in \mathcal{X}^1}\lambda^{x^1x^2}_{y^1y^2}(t),
\quad  x^2,y^2 \in \mathcal{X}^2,\, x^2\ne y^2,\quad
\lambda^2_{x^2x^2}(t)=-\sum_{y^2\in \mathcal{X}^2,y^2\ne x^2}\, \lambda^2_{x^2y^2}(t),\quad
\forall x^2\in \mathcal{X}^2, \ee
then under condition (M2) the functions ${\lambda}^2_{x^2y^2}$ are well defined and locally integrable. It is straightforward to verify that they satisfy \eqref{vn2}, so that weak Markovian consistency holds  with respect to  $X^2$.

 As a matter of fact, it was shown in \cite{BJVV} that conditions (M1) and (M2)  are sufficient for strong Markovian consistency to hold for $X$ with respect to both its components: $X^1$ and $X^2.$ It turns out however, that conditions (M1) and (M2) are too strong; in particular, they are not necessary for strong Markovian consistency to hold for $X$ with respect to its components.

 We now state a theorem providing  sufficient and necessary conditions for strong Markovian consistency of $X$.

\bt\lab{smc} The component $X^1$ of $X$ is a Markov chain with respect to filtration $\mathbb{F}^X$ if and only if
\begin{align}\lab{weak-M1}
\I_{\{X^1_{t}=x^1\!\}}\!\!\!\!\sum_{y^2\in {\mathcal X}^2} \!\!\lambda^{x^1X^2_t}_{y^1y^2}(t)
\!=\! \I_{\{X^1_{t}=x^1\}}\lambda^1_{x^1y^1}(t) \quad {dt \otimes
d\mathbb{P} \text{-a.s.} }  \  \forall x^1\!,y^1\!\in\! {\mathcal
X}^1\!, x^1 \neq y^1
\end{align}
for some locally integrable functions $\lambda^1_{x^1y^1}.$ The
generator of $X^1$ is $\Lambda^1 (t)
=[\lambda^1_{x^1y^1}(t)]_{x^1,y^1\in {\mathcal X}^1}$   with
$\lambda^1_{x^1x^1}$ given by
\[
    \lambda^1_{x^1x^1}(t) = - \sum_{y^1 \in \mathcal{X}^1, y^1 \neq x^1 } \lambda^1_{x^1y^1}(t)
    \qquad
    \forall x^1 \in \mathcal{X}^1
    .
\]
The component $X^2$ of $X$ is a Markov chain with respect to filtration $\mathbb{F}^X$ if and only if
\begin{align}\lab{weak-M2}
\I_{\{X^2_{t}=x^2\!\}}\!\!\!\!\sum_{y^1\in {\mathcal X}^1} \!\!\lambda^{X^1_tx^2}_{y^1y^2}(t)
\!=\! \I_{\{X^2_{t}=x^2\}}\lambda^2_{x^2y^2}(t) \quad {dt \otimes
d\mathbb{P} \text{-a.s.} }  \  \forall x^2\!,y^2\!\in\! {\mathcal
X}^2\!, x^2 \neq y^2
\end{align}
for some locally integrable functions $\lambda^2_{x^2y^2}.$ The
generator of $X^2$ is $\Lambda^2 (t)
=[\lambda^2_{x^2y^2}(t)]_{x^2,y^2\in {\mathcal X}^2}$   with
$\lambda^2_{x^2x^2}$ given by
\[
    \lambda^2_{x^2x^2}(t) = - \sum_{y^2 \in \mathcal{X}^2, y^2 \neq x^2 } \lambda^2_{x^2y^2}(t)
    \qquad
    \forall x^2 \in \mathcal{X}^2
    .
\]
\et
\begin{proof}  We will only give the proof regarding component $X^1$ of $X.$ For the component $X^2$ the proof is analogous.

Assume that  \eqref{weak-M1} holds. Since $X$ is a Markov chain, then, by \eqref{eq:Nu1}, \eqref{nu12} and  \eqref{weak-M1}, for each $x^1$, $y^1 \in \mathcal{X}^1$ the $\mathbb{F}^X$\!\!-compensator of $N^1_{x^1y^1}$ has the form
\[
	\nu\textcolor[rgb]{0.00,0.00,0.00}{^1}_{x^1y^1}(dt) =
\sum_{x^2,y^2\in {\mathcal X}^2} \!\!\lambda^{x^1 x^2}_{y^1y^2}(t)\I_{\{(X^1_{t},X^2_{t})=(x^1,x^2)\!\}}
dt
=
\I_{\{X^1_{t}=x^1\!\}}\!\!\!\!\sum_{y^2\in {\mathcal X}^2} \!\!\lambda^{x^1X^2_t}_{y^1y^2}(t)
dt
=
\I_{\{X^1_{t}=x^1\!\}} \lambda^1_{x^1y^1}(t) dt
\]
for some locally integrable deterministic function $\lambda^1_{x^1y^1}$. Then, by martingale characterization, $X^1$ is a Markov chain with respect to $\mathbb{F}^X$ with generator $\Lambda^1 (t) =[\lambda^1_{x^1y^1}(t)]_{x^1,y^1\in {\mathcal X}^1}.$
Conversely, assume that $X^1$ is a Markov chain with respect to filtration $\mathbb{F}^X$ with generator $\Lambda^1 (t) =[\lambda^1_{x^1y^1}(t)]_{x^1,y^1\in {\mathcal X}^1}$. Then \eqref{weak-M1} follows from martingale characterization, \eqref{eq:Nu1} and \eqref{nu12}.  Indeed, we have
\[
\I_{\{X^1_{t}=x^1\!\}} \lambda^1_{x^1y^1}(t) dt
=
	\nu\textcolor[rgb]{0.00,0.00,0.00}{^1}_{x^1y^1}(dt)
=
\sum_{x^2,y^2\in {\mathcal X}^2} \!\!\lambda^{x^1 x^2}_{y^1y^2}(t)\I_{\{(X^1_{t},X^2_{t})=(x^1,x^2)\!\}}
dt
=
\I_{\{X^1_{t}=x^1\!\}}\!\!\!\!\sum_{y^2\in {\mathcal X}^2} \!\!\lambda^{x^1X^2_t}_{y^1y^2}(t) dt.
\]
\finproof
\end{proof}

\begin{rem}
(i) It is clear that conditions (M1) and (M2) imply conditions \eqref{weak-M1} and \eqref{weak-M2}, respectively. On the other hand, it is clear that conditions \eqref{weak-M1} and \eqref{weak-M2} imply \eqref{vn1} and \eqref{vn2}, respectively.\\
(ii) Even though conditions (M1) and (M2) are stronger that  conditions needed to establish strong Markovian consistency, they are very convenient to use for that purpose. In particular, they can be conveniently used to construct a strong Markov copula (cf. Section \ref{s:smc-1}). In the next section we shall provide operator form of conditions (M1) and (M2).
\end{rem}

\begin{rem}\label{M}
 Ball and Yeo \cite{BY93}  considered time homogeneous
Markov chains with intensity matrix $\Lambda$ satisfying some
additional assumptions (cf. \cite[Condition 2.2]{BY93}). In
\cite[Theorem 3.1]{BY93}, it is proved that the marginal process
$X^1$ of time a homogenous Markov chain $X$ is a time homogenous
Markov  chain in its natural filtration if and only if a condition
equivalent to Condition (M1) holds. However, if we omit the
assumption of time homogeneity, then \cite[Theorem 3.1]{BY93} does
not hold; see our Example \ref{ex2} below. Moreover, assumptions
imposed in \cite{BY93} on $\Lambda$ exclude Markov chains with
absorbing states.
\end{rem}

We shall see in Section \ref{examples} that there exist Markov chains that are weakly Markovian consistent, but not strongly Markovian consistent.

\subsection{Operator interpretation of necessary conditions
for weak Markovian consistency, and of the sufficient condition (M) for strong Markovian consistency}
For $i=1,2$ and $t\geq 0$, we define an operator $Q^i_t$, acting on any function $f$ on $\mathcal X= \mathcal{X}^1 \times \mathcal{X}^2$, by
\be (Q^i_tf)(x^i)=\EP(f(X_t)|X^i_t=x^i),\quad \forall x^i \in {\mathcal X}^i.\ee

We also introduce an extension operator $C^{i,*}$ as follows: for any function $f^i$ on $\mathcal{X}^i$ the function $ C^{i,*} f^i$ is defined on $\mathcal{X} $ by
\[
    (C^{i,*} f^i)(x) = f^{i}(x^i),\quad \forall x=(x^1,x^2) \in {\mathcal X}.
\]

We have the following proposition, which will be important in the next section in the context of weak Markov copulae.

\begin{thm} Fix $i \in \set{1,2}$.   The condition
\be\lab{a1} Q^i_t\Lambda(t)C^{i,*}= \Lambda ^i(t),\ t\geq 0,\ee
where $\Lambda ^i(t)=[{\lambda}^i_{x^iy^i}(t)],$ with functions ${\lambda}^i_{x^iy^i}$ given by \eqref{vn} for $i=1$ and given  by \eqref{vn2} for $i=2$,
 is necessary for weak Markovian consistency with respect to $X^i$.
\end{thm}
\begin{proof} We give the proof for $i=1.$
It is enough to observe that \eqref{vn1} is equivalent to \eqref{a1}.
Indeed, first note that \eqref{a1} is equivalent to the equality
\beqa\label{a1-a}
 (Q^1_t\Lambda(t)C^{1,*} g)(x^1)  = \sum_{y^1 \in \mathcal{X}^1} {\lambda}^1_{x^1y^1}(t) g(y^1)
\eeqa
for an arbitrary function $g$ on $\mathcal{X}^1$ and  $x^1 \in \mathcal{X}^1$.
Now, we rewrite the  left hand side:
\beq
(Q^1_t\Lambda(t)C^{1,*} g)(x^1)
 &=&
 \mathbb{E} \left(  \sum_{(z^1, x^2) \in \mathcal{X}} \I_\set{ X^1_t = z^1, X^2_t = x^2}  \sum_{(y^1, y^2) \in \mathcal{X} } \lambda^{z^1x^2}_{y^1y^2} (t) g(y^1) \bigg| X^1_t = x^1 \right)
 \\
 &=&
  \sum_{ x^2 \in \mathcal{X}^2} \left(\mathbb{E} \left(  \I_\set{ X^2_t = x^2} \big| X^1_t = x^1 \right) \sum_{(y^1, y^2) \in \mathcal{X} } \lambda^{x^1x^2}_{y^1y^2} (t) g(y^1) \right)
  \\
  &
  =&
  \sum_{y^1 \in \mathcal{X}^1 }
  \left(
  \sum_{ x^2  \in \mathcal{X}^2} \sum_{y^2 \in \mathcal{X}^2 } \mathbb{E} \left(  \I_\set{ X^2_t = x^2} \big| X^1_t = x^1 \right) \lambda^{x^1x^2}_{y^1y^2} (t) \right) g(y^1).
\eeq
Since $g$ is arbitrary, \eqref{a1-a} is equivalent to
\[
{\lambda}^1_{x^1y^1}(t)
=
 \sum_{ x^2  \in \mathcal{X}^2} \sum_{y^2 \in \mathcal{X}^2 } \mathbb{E} \left(  \I_\set{ X^2_t = x^2} \big| X^1_t = x^1 \right) \lambda^{x^1x^2}_{y^1y^2} (t),
\]
which is exactly \eqref{vn1}.
\finproof
\end{proof}

In the next two propositions we shall consider an operator interpretation of condition (M) for strong Markovian consistency, and its connection with condition \eqref{a1}.

\begin{pro}
Condition (M1) is equivalent to

\medskip
(N1): There exist generator matrix function $\Lambda^1=[\lambda^1_{x^1y^1}]_{x^1,y^1\in {\mathcal X}^1}$ such that:
\be\lab{lhat1}
C^{1,*}\Lambda^1 (t) =  \Lambda (t) C^{1,*}, \quad  \forall t \geq 0.
\ee

\noindent Condition (M2)  is equivalent to

\medskip

(N2): There exist generator matrix function $\Lambda^2=[\lambda^2_{x^2y^2}]_{x^2,y^2\in {\mathcal X}^2}$ such that:
\be\lab{lhat2.1}
C^{2,*}\Lambda^2 (t) =  \Lambda (t) C^{2,*}, \quad  \forall t \geq 0.
\ee
\end{pro}
\begin{proof} We only prove the first equivalence. The proof of the other one is analogous.
\\
We note that \eqref{lhat1} is equivalent to the equality
\be\lab{lhat1.1}
(C^{1,*}\Lambda^1 (t) g)(x^1,x^2) =  (\Lambda (t) C^{1,*}g)(x^1,x^2),\quad \forall (x^1,x^2)\in \mathcal{X}^1 \times \mathcal{X}^2,
\ee
for an arbitrary function $g$ on $\mathcal{X}^1$. By definition, the right hand side of \eqref{lhat1.1} is
\beq
(\Lambda (t) C^{1,*}g)(x^1,x^2)
&=&
\sum_{(y^1,y^2) \in \mathcal{X} } \lambda^{x^1x^2}_{y^1y^2} (t) (C^{1,*}g)(y^1,y^2) =
\sum_{(y^1,y^2) \in \mathcal{X} } \lambda^{x^1x^2}_{y^1y^2} (t) g(y^1)
\\
&=&
\sum_{y^1 \in \mathcal{X}^1 } \left(\sum_{y^2 \in \mathcal{X}^2 }  \lambda^{x^1x^2}_{y^1y^2} (t)\right) g(y^1),
\eeq
and the left hand side of \eqref{lhat1.1} is given by
\beq
    (C^{1,*}\Lambda^1 (t) g)(x^1,x^2)
    &=&
    \sum_{y^1 \in \mathcal{X}^1 }
    \lambda^1_{x^1y^1} (t) g(y^1).
\eeq
Since $g$ is arbitrary, we obtain that (N1) is equivalent to existence of matrix function $\Lambda^1=[\lambda^1_{x_1y_1}]_{x^1,y^1\in {\mathcal X}^1}$ such that for each $ t \geq 0$ we have
\be\label{eq:M3}
    \lambda^1_{x^1y^1} (t) = \sum_{ y^2 \in \mathcal{X}^2} \lambda^{x^1x^2}_{y^1 y^2} (t)
    \qquad
    \forall x^1, y^1 \in \mathcal{ X   }^1,\
    \forall x^2 \in \mathcal{ X   }^2.
\ee

Hence, using the fact that $\Lambda$ is  the generator of a Markov chain we  see that
(N1) is equivalent to (M1). Finally, note also that in a view of \eqref{eq:M3} it is straightforward to verify that matrix function $\Lambda^1$ is a valid generator matrix.
\finproof
\end{proof}
\begin{pro}\label{pro}
Condition \eqref{lhat1} implies \eqref{a1} for $i=1$ and condition \eqref{lhat2.1} implies \eqref{a1} for $i=2$.
\end{pro}
\begin{proof}
Since  $Q^i_tC^{i,*}={\textrm {Id}}$ for $i=1,2$, we have
\beq Q^i_t\Lambda(t)C^{i,*}=Q^i_tC^{i,*} \Lambda^i (t)= \Lambda^i (t),\ t\geq 0, \quad i=1,2. \eeq \finproof
\end{proof}
\begin{rem}Another possible proof of Proposition \ref{pro} is the following: Conditions \eqref{lhat1} and \eqref{lhat2.1} 
are sufficient for strong Markovian consistency of $X($see Remark \ref{M}$)$, which implies weak Markovian consistency of $X$, for which \eqref{a1} is a necessary condition.
\end{rem}
\begin{rem}
In the case of time homogeneous Markov processes, conditions analogous to \eqref{lhat1} and \eqref{lhat2.1} have been previously studied in \cite{BVV} and \cite{Andrea}, and it has been shown that they are sufficient for strong Markovian consistency.
So, \eqref{lhat1} and \eqref{lhat2.1} imply that each coordinate of the Markov process in question is a Markov process with respect to $\mathbb{F}^X$. It is worth noting that \eqref{lhat1} and \eqref{lhat2.1} agree with (10.60) of Dynkin \cite{Dyn}, if the latter is applied to $f$  being a component projection function.
\end{rem}

\begin{rem} The operator conditions \eqref{lhat1} and \eqref{lhat2.1} for strong Markovian consistency can be interpreted in the context of martingale characterization of Markov chains.

Let $C^{i},$ $i=1,2,$ be the projection from $\mathcal{X}^1 \times \mathcal{X}^2$ on the $i$th component.
Fix $i \in \set{1,2}$ and $0\leq s\leq t$. Since $X$ is a Markov chain, for any function $f^i$ on ${\mathcal X}^i$ we have the representation
\be
C^{i,*}f^i(X_t) = C^{i,*}f^i(X_s) + \int _s ^t (\Lambda(u)(C^{i,*}f^i))(X_u)du + M^{C^{i,*},f^i}_t-M^{C^{i,*},f^i}_s,
\ee
where $M^{C^{i,*},f^i}$ is a martingale  with respect to  $\ff^X$. Thus,
\be\lab{A}
f^i(C^iX_t) = f^i(C^iX_s) + \int _s ^t (\Lambda(u)(C^{i,*}f^i))(X_u)du +  M^{C^{i,*},f^i}_t-M^{C^{i,*},f^i}_s.
\ee
If conditions \eqref{lhat1} and \eqref{lhat2.1} hold then we may rewrite \eqref{A} as
\be
\lab{Ai} f^i(X^i_t) = f^i(X^i_s) + \int _s ^t (\Lambda^i(u)f^i)(X^i_u)du + M^{C^{i,*},f^i}_t-M^{C^{i,*},f^i}_s,
\ee
which shows that $X^i$ is a Markov chain  with respect to  $\ff^X.$
\end{rem}

\subsection{When Does Weak Markov Consistency Imply Strong Markov Consistency?}

It is well known that if a process $X$ is a $\P$-Markov chain  with respect to  a filtration $\ff$,
and if it is adapted with respect to a filtration $\hat \ff \subset \ff,$ then $X$ is a
$\P$-Markov chain  with respect to  $\hat \ff.$ However, the converse is not true in general.
Nevertheless, if $X$ is a $\P$-Markov chain with respect to $\hat \ff$, and $\hat \ff$ is $\P$-{\it immersed}
 in $\ff$ \footnote{We say that a filtration $\hat \ff$ is $\P$-{\it immersed} in a filtration $\ff$ if $\hat \ff \subset \ff$ and every $(\P,\hat \ff)$-local-martingale is a  $(\P,\ff)$-local-martingale. 
}, then we can deduce from the martingale characterization of Markov chains that $X$ is also a  $\P$-Markov chain with respect to $\ff$.

Thus, if  $\ff^{X^i}$ is $\P$-immersed in $\ff^X$, then weak Markovian consistency of $X$  with respect to  $X^i$ will imply strong Markovian consistency of $X$  with respect to  $X^i.$ In the following theorem we demonstrate that in fact this property is equivalent to $\P$-immersion between $\ff^{X^i}$ and $\ff^X$, given that weak Markovian consistency holds.

\begin{thm}
    Assume that $X$ satisfies the weak Markovian consistency condition with respect to $X^i$.   Then $X$ satisfies the strong Markovian consistency condition if and only if $\mathbb{F}^{X^i}$ is $\P$-{\it immersed} in $\ff^X$.
\end{thm}
\begin{proof}
    $"\Longrightarrow"$
    We give a proof in the case of $i = 1$.
        By Proposition \ref{prop:comp} the process
    \[
        M^1_{x^1y^1}(t) :=  N^1_{x^1y^1}(t) - \int_{(0,t]}  \wh{\nu}^1_{x^1y^1} (ds)
    \]
    is an $\ff^{X^1}$-martingale for every $x^1 \neq y^1$ since $X^1$ is a Markov process  with respect to  its own filtration.
    By  Jeanblanc, Yor and Chesney \cite[Proposition 5.9.1.1]{JYCh2009} it is sufficient to show that every $\ff^{X^1}$-square integrable martingale $Z$ is  also an $\ff^X$-square integrable martingale under $\P$.
    Using the martingale representation theorem (see Rogers and Williams \cite[Theorem 21.15]{RW2000}) we have
    \begin{align}\label{eq:mart-rep}
        Z_t =  Z_0 + \sum_{ x^1 \neq y^1} \int_{(0,t]} g(s,x^1,y^1, \omega ) ( N^1_{x^1y^1}(ds) - \wh{\nu}^1_{x^1y^1} (ds)  )
    \end{align}
    for some function $g: (0,\infty) \times \mathcal{X}^1 \times  \mathcal{X}^1 \times \Omega \rightarrow \mathbb{R},$ such that for every $x^1, y^1$ the mapping $(t,\omega) \mapsto g(t,x^1,y^1, \omega)$ is $\ff^{X^1}$-predictable and $g(t,x^1,x^1,\omega) = 0, $  $\P$-a.s.$\, .$ The $\ff^{X^1}$-oblique bracket of $M^1_{x^1y^1}$ (i.e.  the $\ff^{X^1}$\!\!\!\!-compensator of $(M^1_{x^1,y^1})^2$) is equal to  $(\int_0^t \wh{\nu}^1_{x^1y^1} (ds))_{t \geq 0 }$, and therefore $g$ satisfies the integrability condition
    \begin{align}\label{eq:iso-cond}
        \mathbb{E} \left( \sum_{ x^1 \neq y^1} \int_{(0,T]} |g(s,x^1,y^1 )|^2  \wh{\nu}^1_{x^1y^1} (ds)  \right) < \infty
        \qquad \forall \ T>0.
    \end{align}
    From the assumption that weak Markovian consistency implies strong Markovian consistency  we infer that
    $X^1$ is a Markov chain with respect to $\ff^X$, and therefore $M^1_{x^1y^1}$ are $\ff^X$-martingales for every $x^1 \neq y^1$. Moreover, the $\ff^{X}$-oblique bracket of $M^1_{x^1y^1}$ is also equal to $(\int_0^t \wh{\nu}^1_{x^1y^1} (ds) )_{t \geq 0 }$, and obviously for every $x^1, y^1$ the mapping $(t,\omega) \rightarrow g(t,x^1,y^1, \omega)$ is $\ff^{X}$-predictable.
    Hence using \eqref{eq:mart-rep}  and \eqref{eq:iso-cond} we deduce that $Z$ is also an $\ff^X$-square integrable martingale.

    \noindent $"\Longleftarrow"$ Assume that $\ff^{X^i}$ is immersed in  $\ff^X$. Weak Markovian consistency for $X^1$ implies that the process $M^1_{x^1y^1}$ is an $\ff^{X^1}$-martingale for every $x^1 \neq y^1$. By immersion we know that $M^1_{x^1y^1}$ are $\ff^{X}$-martingales for every $x^1 \neq y^1$ and therefore Proposition \ref{prop:comp} implies that $X^1$ is a Markov process with respect to $\ff^X$.
\finproof
\end{proof}

\section{Markov copulae}\lab{s:wmc}

We now turn to the problem of constructing a multivariate finite Markov chain whose components are finite univariate Markov chains with given generator matrices.

This problem was previously studied in \cite{BJVV} and \cite{BJN2}, for example, in the context of strong Markovian consistency. This meant that the components of the multivariate Markov chain constructed were Markovian  both in  their own filtrations and in the filtration of the entire chain. Thus, essentially, these references dealt with constructing of what we shall term here  \emph{strong Markov copulae}.

In this paper, we shall additionally be concerned with \emph{weak
Markov copulae} in the context of finite Markov chains. It will be
seen that any strong Markov copula is also a weak Markov copula.

As in the previous section, in order to simplify the notation we shall consider bivariate processes $X$ only.

\subsection{Strong Markov copulae}\lab{s:smc-1}

The key observation leading to the concept of strong Markov copula is the following:  Let there be given two generator functions  $\Lambda^1(t)=[\lambda^1_{x^1y^1}(t)]$ and $\Lambda^2(t)=[\lambda^2_{x^2y^2}(t)]$, and suppose that there exists a valid generator  matrix function $\Lambda(t)=[\lambda^{x^1x^2}_{y^1y^2}(t)]_{x^1,y^1\in \mathcal{X}^1, x^2,y^2\in
\mathcal{X}^2}$ satisfying \eqref{lhat} for every $x^2 \in \mathcal{X}^2$, and satisfying \eqref{lhat2} for every $x^1 \in \mathcal{X}^1$. Then, Condition {\sf (M)} is clearly satisfied, so that (cf. Remark \ref{M})
strong Markovian consistency holds for the Markov chain, $X$ generated by $ \Lambda(t)$.

Note that, typically, system \eqref{lhat} and \eqref{lhat2}, considered as a system with given $\Lambda^1(t)=[\lambda^1_{x^1y^1}(t)]$ and $\Lambda^2(t)=[\lambda^2_{x^2y^2}(t)]$ and with unknown $\Lambda(t)=[\lambda^{x^1x^2}_{y^1y^2}(t)]_{x^1,y^1\in \mathcal{X}^1, x^2,y^2\in
\mathcal{X}^2}$, contains many more unknowns (i.e., $\lambda^{x^1x^2}_{y^1y^2}(t)$, $x^1,y^1\in \mathcal{X}^1, x^2,y^2\in
\mathcal{X}^2$) than it contains equations. In fact, given that the cardinalities of $\mathcal{X}^1$ and
$\mathcal{X}^2$ are $K_1$ and $K_2,\, $
respectively, the system consists of $\KS(\KS-1)+\KO(\KO-1)$
equations in $\KS\KO(\KS\KO-1)$ unknowns.

Thus, in principle, one can create several bivariate Markov chains
$X$ with margins $X^1$ and $X^2$ that are  Markovian in the filtration of $X$, and such that the law of $X^i$ agrees with the law of a given Markov chain $Y^i$, $i=1,2.$ Thus, indeed, the system \eqref{lhat} and \eqref{lhat2} essentially serves as a "copula"\footnote{We use the term "copula" in analogy to classical copulae for probability distributions of finite-dimensional random variables (cf. e.g. \cite{Nelsen}). See also discussion in Section \ref{sklar}.} between
the Markovian margins $Y^1$, $Y^2$ and the bivariate Markov chain $X.$
This observation leads to the following definition,

\bd Let $Y^1$ and $Y^2$ be two Markov chains with values in $\mathcal{X}^1$ and $\mathcal{X}^2$, and with generators $\Lambda^1(t)=[\lambda^1_{x^1y^1}(t)]$ and $\Lambda^2(t)=[\lambda^2_{x^2y^2}(t)]$. A {\sf
Strong Markov Copula} between the Markov chains $Y^1$ and $Y^2$ is any
solution to \eqref{lhat} and \eqref{lhat2} such that the matrix
function $\Lambda(t)=[\lambda^{x^1x^2}_{y^1y^2}(t)]_{x^1,y^1\in \mathcal{X}^1, x^2,y^2\in
\mathcal{X}^2}$, with $\lambda^{x^1x^2}_{x^1x^2}(t)$ given as
\be \lambda^{x^1x^2}_{x^1x^2}(t)=- \sum_{(z^1,z^2)\in \mathcal{X}^1\times \mathcal{X}^2,\, (z^1,z^2)\ne (x^1,x^2)} \lambda^{x^1x^2}_{z^1z^2}(t),\ee
correctly defines the infinitesimal generator function of a Markov
chain with values in $\mathcal{X}^1\times \mathcal{X}^2.\, $  \ed

Thus, any strong Markov copula between Markov chains $Y^1$ and $Y^2$ produces a bivariate Markov chain, say $X=(X^1,X^2)$, such that
\begin{itemize}
\item the components $X^1$ and $X^2$ are Markovian in the filtration of $X$,
\item the transition laws of $X^i$ is the same as the transition laws of $Y^i$, $i=1,2$,
\item If, in addition, the initial law of $X^i$ is same as the initial laws of $Y^i$,  then, the law  $X^i$ is the same as the law of $Y^i$, $i=1,2.$  In this case, according to terminology of \cite{BJN}, the process $X$ satisfies the strong Markovian consistency condition relative to $Y^1$ and $Y^2$.
\end{itemize}

It is clear that there exists at least one solution to
\eqref{lhat} and \eqref{lhat2} such that the matrix function
$\Lambda(t)=[\lambda^{x^1x^2}_{y^1y^2}(t)]_{x^1,y^1\in
\mathcal{X}^1, x^2,y^2\in \mathcal{X}^2}$ is a valid generator
matrix. This solution correspond to the case of independent
processes $X^1$ and $X^2$. In this case we have $\Lambda(t)=I^1
\hat{\otimes} \Lambda^2(t) + \Lambda^1(t) \hat{\otimes} I^2$ where
$A \hat{\otimes} B$  denotes tensor product of operators $A$ and
$B$ (see Ryan \cite{Ryan2000}), and where $I^i$ is identity
operator on $\mathcal{X}^i$. Matrix $\Lambda(t) $ that corresponds
to two independent processes can be also written more explicitly
\[
    \lambda^{x^1x^2}_{y^1y^2}(t) =
    \left\{
      \begin{array}{ll}
        \lambda^1_{x^1x^1}(t) + \lambda^2_{x^2x^2}(t), & y^1 = x^1, y^2 = x^2, \\
        \lambda^1_{x^1y^1}(t), & y^1 \neq x^1, y^2 = x^2, \\
        \lambda^2_{x^2y^2}(t), & y^2 \neq x^2, y^1 = x^1, \\
        0, & \hbox{otherwise.}
      \end{array}
    \right.
\]

\subsection{Weak Markov Copulae}

The concept of weak Markov copula corresponds to the concept of weak Markovian consistency. We do not have any clear analytical characterization of the latter property, analogous to condition (M) that is sufficient for strong Markovian consistency.

Consequently, the concept of weak Markov copula is much more intricate than that of strong Markov copula, because  it involves both probabilistic and analytical (indeed, algebraic in our case) characterizations.

\bd  Let $Y^1$ and $Y^2$ be two Markov chains with values in $\mathcal{X}^1$ and $\mathcal{X}^2$, and with generators $\Lambda^1(t)=[\lambda^1_{x^1y^1}(t)]$ and $\Lambda^2(t)=[\lambda^2_{x^2y^2}(t)]$, respectively. A {\sf
Weak Markov Copula} between $Y^1$ and $Y^2$ is any
matrix function $\Lambda(t)=[\lambda^{x^1x^2}_{y^1y^2}(t)]_{x^1,y^1\in \mathcal{X}^1, x^2,y^2\in
\mathcal{X}^2}$ 
that satisfies the following conditions:
\begin{description}
\item[\bf (WMC1)] $\Lambda(t)$  correctly properly defines the  infinitesimal generator  of a bivariate Markov
chain, say $X=(X^1,X^2)$, with values in $\mathcal{X}^1\times \mathcal{X}^2\, $,
\item[\bf (WMC2)] Conditions \eqref{vn} and \eqref{vn2} are satisfied, so that $X$ is weakly Markovian consistent.
\end{description}
\ed

Thus, any weak Markov copula between the Markov chains $Y^1$ and $Y^2$ produces a bivariate Markov chain, say $X=(X^1,X^2)$, such that
\begin{itemize}
\item the components $X^1$ and $X^2$ are Markovian in their own filtrations, but not necessarily Markovian in the filtration of $X$, and

\item the transition laws of $X^i$ is the same as the transition laws of $Y^i$, $i=1,2$,
\item If, in addition, the initial law of $X^i$ is same as the initial laws of $Y^i$,  then, the law $X^i$ is the same as the law of $Y^i$, $i=1,2.$  In this case, we say that the process $X$ satisfies the weak Markovian consistency condition relative to $Y^1$ and $Y^2$.
\end{itemize}

It is clear that any strong Markov copula between $Y^1$ and $Y^2$ is also a weak Markov copula between $Y^1$ and $Y^2$.

 A possible way of constructing a \textbf{weak-only} Markov copula, that is a  weak Markov copula, \textbf{which is not} a strong Markov copula, is to start with the necessary condition \eqref{a1} and to find a generator matrix   $\Lambda(t)$ that satisfies this condition with given  $\Lambda^1(t)$ and  $\Lambda^2(t)$. Typically,  the matrix $\Lambda(t)$ found will generate a Markov chain satisfying the weak Markovian consistency condition relative to the Markov chains $Y^1$ and $Y^2$ generated by  $\Lambda^1(t)$ and  $\Lambda^2(t)$, respectively. This approach will be illustrated in Example \ref{ex2} below.

\begin{rem}\label{fin}
It needs to be strongly stressed that the issue of constructing {weak-only} Markov copulae is very important from the practical point of view. For example, it is important
in the context of credit risk management since weak Markov copulae allow for modeling of default contagion between individual obligors and the rest of the credit pool (cf. \cite{BCCH} for a discussion); this kind of contagion is precluded in the context of strong Markov copulae. Thus, weak Markov copulae make it possible to tackle two critical modeling requirements:
\begin{itemize}
\item
They make it possible to model contagion between credit events in credit portfolios; equally importantly, they allow for modeling contagion between failure events in complex manufacturing systems;
\item They make it possible to separate calibration of the model to univariate data (credit default spreads, for example), from calibration of the model to multivariate data (spreads on credit portfolio contracts, such as collateralized loan obligations or collateralized debt obligations). This aspect of the Markov copula theory is of fundamental importance for efficient calibration of a model to market data. In \cite{BVV1} and \cite{BCCH} (see also references therein), the strong Markov copula theory was successfully applied to separate calibration of dependence in the pool of 125 obligors (constituting an iTraxx index), from the calibration of univariate characteristics of the individual obligors. We are currently working on using the weak-only Markov copulae for such purpose.
\end{itemize}
\end{rem}


\subsection{Classical copulae theory vs Markov copulae theory}\lab{sklar}

It is useful to relate the concept of Markov copulae to the classical concept of \emph{copula function} used in probability to construct multivariate, finite dimensional random variables, with given marginal distributions.

Recall that a function $C: [0,1]^N \rightarrow [0,1]$ is an $N-$copula if, and only if, the following properties hold:

\begin{enumerate}
	 \item for every $j \in \set{1,2, . . . ,N}$, $C(1,\ldots,1,u_j,1,\ldots,1) = u_j$; 

  \item $C$ is isotonic, that is $C(u) \leq C(v)$ for all $u,v \in  [0,1]^N, u \leq v$;

  \item $C$ is N-increasing, that is
    \[
        \sum_{w \in \set{u_1, v_1}\times \ldots \times \set{ u_N, v_N}} (-1)^{ \# \set{j: w_j =v_j} } C(w) \geq 0
    \]
for all $u, v \in [0,1]^N$, $u \leq v$.
\end{enumerate}

Let now $U_1,\ldots,U_N$ be real valued random variables, with the corresponding cumulative distribution functions $F_1, F_2, \ldots , F_N,$ and let  $C$ be an $N$-copula. Next, let the function $F : \mathbb{R}^N \rightarrow [0,1]$ be defined by
\be\label{cop}
F(u_1,u_2, . . . ,u_N) =C(F_1(u_1),F_2(u_2), \ldots ,F_N(u_N)).
\ee
It is the classical result due to Sklar \cite{Sklar1973} that $F$ is a cumulative distribution function of an $\mathbb{R}^N$-valued random variable, say $W=(W_1,\ldots,W_N)$, such that the law of $W_n$ is the same as the law of $U_n$, $n=1,2,\ldots,N.$ In other words, $F$ is an $N$-variate distribution function with margins $F_1, F_2, \ldots , F_N$.

Now, we have the following  analogies between the classical copula theory and the Markov copulae (below, we use our convention that $N=2$):
 \begin{itemize}

 \item random variables $U_1$ and $U_2$ correspond to Markov chains $Y^1$ and $Y^2$, random variables $W_1$ and $W_2$ correspond to Markov chains $X^1$ and $X^2$, and random variable $W=(W_1,W_2)$ corresponds to Markov chain $X=(X^1,X^2),$

 \item the generator functions $\Lambda^1(\cdot)=[\lambda^1_{x^1y^1}(\cdot)]$ and $\Lambda^2(\cdot)=[\lambda^2_{x^2y^2}(\cdot)]$, showing in equations \eqref{lhat} and \eqref{lhat2} and in equations \eqref{vn} and \eqref{vn2}, are analogous to the marginal distributions $F_1$ and $F_2$ showing in \eqref{cop},

 \item in the case of strong Markov copula, the equations \eqref{lhat} and \eqref{lhat2} and any of their solutions, say $\Lambda(\cdot)=[\lambda^{x^1x^2}_{y^1y^2}(\cdot)]_{x^1,y^1\in \mathcal{X}^1, x^2,y^2\in
\mathcal{X}^2}$, which produces a valid Markov chain, is analogous to the pair $(F,C)$ in \eqref{cop},

 \item in the case of weak Markov copula, the equations \eqref{vn} and \eqref{vn2} and any of their solutions, say $\Lambda(\cdot)=[\lambda^{x^1x^2}_{y^1y^2}(\cdot)]_{x^1,y^1\in \mathcal{X}^1, x^2,y^2\in
\mathcal{X}^2}$, which produces a valid Markov chain, is analogous to the pair $(F,C)$ in \eqref{cop}.\,

\end{itemize}

It needs to be stressed that we use the term "copula", in \textsl{Markov copula}, because of the above correspondences, and, really, for reason of tradition. In general, there is no \emph{copula functional} that would map marginal Markov processes $X^n$ to a multivariate Markov process $X$ (cf. discussion of this issue given in \cite{BJVV}).

\section{Examples}\label{examples}
As before, we take $N=2$ in the examples below. We shall present
examples illustrating
\begin{itemize}
\item Construction of a strong Markov copula (Example \ref{ex1}),
i.e., a construction of a two dimensional Markov chain $X=(X^1,X^2)$
with components $X^1$ and  $X^2$ that are Markovian in the filtration
of $X$, and such that the transition laws of $X^i$ agree with the transition laws of a
given Markov chain $Y^i$, $i =1,2$.
\item Construction of a weak-only Markov copula (Example \ref{ex2}),
i.e.,  a construction of a two dimensional Markov chain
$X=(X^1,X^2)$ with the components $X^1$ and  $X^2$ that are Markovian in
their own filtrations, but are not Markovian in the filtration of $X$,
and such that the transition laws of $X^i$ agree with the transition laws of a given
Markov chain $Y^i$, $i =1,2$.
\item Existence of a Markov chain for which weak Markovian consistency does not hold, that is, a Markov chain that can't serve as a weak Markov copula (Example \ref{ex3}).
In this example, component $X^2$ of Markov chain $X=(X^1,X^2)$ is shown to be not Markovian in its own filtration.
\end{itemize}
 \bex\lab{ex1} Let us consider two processes, $Y^1$ and $Y^2$,
that are time-homogeneous Markov chains, each taking values in the
state space $\{0,1\},$ with respective generators
\begin{eqnarray}\label{matrgen-1}
\Lambda^1=\bordermatrix{&0&1\cr
                0&-(a+c) &  a+c \cr
                1& 0  &   0 }
\end{eqnarray}
and
\begin{eqnarray}\label{matrgen-2}
\Lambda^2=\bordermatrix{&0&1\cr
                0&-(b+c) &  b+c \cr
                1& 0  &   0 },
\end{eqnarray}
for $a,b,c\geq 0.$

We shall first consider the system of equations \eqref{lhat1} and \eqref{lhat2.1} for this example. In this case we identify $C^{i,*},\, i=1,2,$ with the matrices
\be
 C^{1,*} = \left( \begin{array}{cc} 1 & 0 \\ 1 & 0 \\ 0 & 1 \\ 0 & 1 \end{array} \right)\quad  \textrm{and}\quad  C^{2,*} = \left( \begin{array}{cc} 1 & 0 \\ 0 & 1 \\ 1 & 0 \\ 0 & 1 \end{array} \right).
\ee
It can be easily checked that the matrix $\Lambda $ below satisfies \eqref{lhat1} and \eqref{lhat2.1}:
\be\lab{L} \Lambda=\bordermatrix{&(0,0)&(0,1)&(1,0)&(1,1)\cr
                (0,0)&-(a+b+c) &  b  & a & c\cr
                (0,1)& 0  &  -(a+c) & 0 & a+c\cr
                (1,0)& 0 & 0 & -(b+c) & b+c\cr
                (1,1)& 0  &   0 &0 & 0}.\ee

Thus, according to the theory of Section  \ref{s:wmc}, $\Lambda$ is a strong Markov copula between $Y^1$ and $Y^2$. Nevertheless, it will be instructive to verify this directly. Towards this end, let us consider the bivariate Markov chain $X=(X^1,X^2)$ on the state space
\begin{equation*}
E=\left\{ {\left( 0,0\right) },{\left(
0,1\right) },{\left( 1,0\right) },{\left( 1,1\right) }\right\}
\end{equation*}
generated by the matrix $\Lambda$ given by \eqref{L}.
We first compute the transition probability matrix for $X$, for $t\geq 0$:
\begin{eqnarray*}
P(t)\!=\!\!\left(
\begin{array}{cccccccc}
  e^{-(a+b+c)t}   & e^{-(a+c)t}(1\!-\!e^{-bt})  & e^{-(b+c)t}(1\!-\!e^{-at})  &
  e^{-(a+b+c)t}\!-\!e^{-(b+c)t}\!-\!e^{-(a+c)t}\!+\!1  \\[2mm]
  0 & e^{-(a+c)t}  & 0 & 1-e^{-(a+c)t}   \\[2mm]
  0 & 0 & e^{-(b+c)t} &  1-e^{-(b+c)t}   \\[2mm]
  0 & 0 & 0 & 1
\end{array}\right)
\end{eqnarray*}
Thus, for any $t \geq 0$,
\allowdisplaybreaks
\begin{align*}
& \lim_{h\rightarrow 0}\frac{P(X^2_{t+h}=0|X^2_t =0)-1}{h}
=-(b+c).
\end{align*}
Similarly, for any $t \geq 0$,
\[
\lim_{h\rightarrow 0}\frac{P(X^1_{t+h}=0|X^1_t =0)-1}{h} =-(a+c).
\]
It is clear that $X^1$ and $X^2$ are Markov chains in their own filtrations (as both chains are absorbed in state $1$). From the above calculations we see that the generator of $X^i$ is $\Lambda ^i$, $i=1,2.$

To verify that $\Lambda$ is a strong Markov copula between $Y^1$ and $Y^2$, it remains to show that components $X^1$ and $X^2$ are Markovian in the filtration of $X$. This can also be verified by direct computations: indeed,
\begin{align*}
  \lim_{h\rightarrow 0}\frac{P(X^1_{t+h}=0|X^1_t =0,X^2_t=0)-1}{h}&=\lim_{h\rightarrow 0}\frac{P(X^1_{t+h}=0|X^1_t =0,X^2_t=1)-1}{h}
 \\
 &=-(a+c) =\lim_{h\rightarrow 0}\frac{P(X^1_{t+h}=0|X^1_t =0)-1}{h} ,
\end{align*}
or, equivalently,
$$P(X^1_{t+h}=0|X^1_t =0,X^2_t=0)=P(X^1_{t+h}=0|X^1_t =0,X^2_t=1)=P(X^1_{t+h}=0|X^1_t =0)=e^{-(a+c)h},
$$
so that condition \eqref{eq:markov-cons}  is satisfied for $X^1$, and similarly for $X^2.$

Note that in accordance with the concept of strong Markovian consistency, the transition intensities and transition probabilities for $X^1$ do not depend on the state of $X^2$:
\begin{itemize}
\item No matter what the state of $X^2$ is, whether $0$ or $1$, the intensity of transition of $X^1$ from $0$ to $1$ is equal to $a+c.$
\item The  transition  probability of $X^1$ from $0$ to $1$ in $t$ units of time, no matter what the state of $X^2$ is, is equal to
       \[
     e^{-(b+c)t}(1-e^{-at})+e^{-(a+b+c)t}-e^{-(b+c)t}-e^{-(a+c)t}+1 = 1-e^{-(a+c)t}.
    \]
\end{itemize}
An analogous observation holds for $X^2$.
Finally, note that $Y^1$ and $Y^2$ are independent if and only if $c=0$.
\eex

\bex\lab{ex2} Let us consider two processes, $Y^1$ and $Y^2$, that are  Markov chains, each taking values in the state space $\{0,1\},$ with respective generator functions
\begin{eqnarray*}
\Lambda^1(t)=\left(
\begin{array}{cccccccc}
  -(a+c)+\alpha(t)   & a+c-\alpha(t)    \\[2mm]
   0 & 0
\end{array}\right)
\end{eqnarray*}
and
\begin{eqnarray*}
\Lambda^2(t)=\left(
\begin{array}{cccccccc}
  -(b+c)+\beta(t)   & b+c-\beta(t)   \\[2mm]
   0 & 0
\end{array}\right),
\end{eqnarray*}
where
\[ \alpha(t)= c\cdot \frac{ e^{-at}(1-e^{-(b+c)t})\frac{b}{b+c} }{e^{-(a+b+c)t}  + e^{-at}(1-e^{-(b+c)t})\frac{b}{b+c}},\quad \beta(t)=
c\cdot \frac{e^{-bt}(1-e^{-(a+c)t})\frac{a}{a+c} }{ e^{-(a+b+c)t} + e^{-bt}(1-e^{-(a+c)t})\frac{a}{a+c}  },
\]
for $a,b,c\geq 0.$

Here we shall seek a weak Markov copula for $Y^1$ and $Y^2$.  Thus  we shall investigate the necessary condition  \eqref{a1}. Towards this end we first note that in this example the matrix representation of the operator $Q^1_t$ takes the form
\begin{eqnarray*}\label{Q1S}
\tiny
Q^1_t\!=\!\!\left(
\begin{array}{cccccccc}
 P(X^1_t\!=\!0,X^2_t\!=\!0|X^1_t\!=\!0)  & P(X^1_t\!=\!0,X^2_t\!=\!1|X^1_t\!=\!0)  &  P(X^1_t\!=\!1,X^2_t\!=\!0|X^1_t\!=\!0)  &
   P(X^1_t\!=\!1,X^2_t\!=\!1|X^1_t\!=\!0) \\[2mm]
   P(X^1_t\!=\!0,X^2_t\!=\!0|X^1_t\!=\!1) &  P(X^1_t\!=\!0,X^2_t\!=\!1|X^1_t\!=\!1)&   P(X^1_t\!=\!1,X^2_t\!=\!0|X^1_t\!=\!1) &   P(X^1_t\!=\!1,X^2_t\!=\!1|X^1_t\!=\!1)  \end{array}\right),
\end{eqnarray*}
and similarly for $Q^2_t.$ It turns out that a solution to the necessary condition  \eqref{a1} is a valid generator matrix
\begin{eqnarray}\label{matrgenL}
\Lambda=\left(
\begin{array}{cccccccc}
  -(a+b+c)  & b  & a   &
  c  \\[2mm]
  0 & -a  & 0 & a   \\[2mm]
  0 & 0 & -b & b  \\[2mm]
  0 & 0 & 0 & 0
\end{array}\right),
\end{eqnarray}
where $a,b\geq 0$ and $ c> 0$. Verification of this is straightforward, but computationally intensive, and can be obtained from the authors on request. 

Since condition \eqref{a1} is just a necessary condition for weak Markovian consistency, the matrix $\Lambda $ in \eqref{matrgenL} may not be a weak Markov copula for  $Y^1$ and $Y^2$. This  has to be verified by direct inspection.

Let us consider the bivariate Markov chain $X=(X^1,X^2)$ on the state space
\begin{equation*}
E=\left\{ {\left( 0,0\right) },{\left(
0,1\right) },{\left( 1,0\right) },{\left( 1,1\right) }\right\}
\end{equation*}
generated by the matrix $\Lambda$ given by \eqref{matrgenL}.

Arguing as in the previous example,  it is clear that the components $X^1$ and $X^2$ are Markovian in their own filtrations. We shall show that:
\begin{itemize}
\item  $X^1$ and $X^2$ are NOT Markovian in the filtration $\ff^X,$ and

\item the generators of $X^1$ and $X^2$ are given by \eqref{matrgen-3} and \eqref{matrgen-4}, respectively.

\end{itemize}


We first compute the transition probability matrix for $X$, for $t\geq 0$:
\begin{eqnarray*}
P(t)\!=\!\!\left(\!\!
\begin{array}{cccccccc}
  e^{-(a+b+c)t}  & e^{-at}(1\!-\!e^{-(b+c)t})\frac{b}{b+c} & e^{-bt}(1\!-\!e^{-(a+c)t})\frac{a}{a+c}   & 1\!+\!
  e^{-(a+b+c)t}(\frac{a}{a+c}\!-\!\frac{c}{b+c})\!-\!\frac{a}{a+c}e^{-bt}\!-\!\frac{b}{b+c}e^{-at}  \\[2mm]
 0 &   e^{-at}& 0 &  1-e^{-at}   \\[2mm]
  0 & 0 & e^{-bt}  & 1-e^{-bt}   \\[2mm]
   0 & 0 & 0 & 1
\end{array}\!\!\right).
\end{eqnarray*}
It follows that
\begin{align*}
P(X^1_{t+h}=0|X^1_t =0,X^2_t=0)&=e^{-(a+b+c)t}+ e^{-at}(1-e^{-(b+c)t})\frac{b}{b+c}
\\
& \ne P(X^1_{t+h}=0|X^1_t =0,X^2_t=1)=e^{-at}
\end{align*}
unless $c=0,$ which is the case of independent $X^1$ and $X^2.$ Thus, in general, $X^1$ is NOT a Markov process in the full filtration. Similarly for $X^2.$

We shall now compute the generator function for $X^2.$ As in the previous example, for any $t \geq 0$,\allowdisplaybreaks
\begin{align*}
&\lim_{h\rightarrow 0}\frac{P(X^2_{t+h}=0|X^2_t =0)-1}{h}
=-(b+c)+c\frac{P(X^1_t=1,X^2_t=0)}{P(X^2_t=0)}.
\end{align*}
Similarly, for any $t \geq 0$,
\[
\lim_{h\rightarrow 0}\frac{P(X^1_{t+h}=0|X^1_t =0)-1}{h} =-(a+c)+c\frac{P(X^1_t=0,X^2_t=1)}{P(X^1_t=0)}.
\]
Thus, both $X^1$ and $X^2$ are time-inhomogeneous Markov chains with generator functions, respectively,
\begin{eqnarray}\label{matrgen-3}
A^1(t)=\left(
\begin{array}{cccccccc}
  -(a+c)+c\frac{P(X^1_t=0,X^2_t=1)}{P(X^1_t=0)}   & a+c-c\frac{P(X^1_t=0,X^2_t=1)}{P(X^1_t=0)}   \\[2mm]
   0 & 0
\end{array}\right)
\end{eqnarray}
and
\begin{eqnarray}\label{matrgen-4}
A^2(t)=\left(
\begin{array}{cccccccc}
  -(b+c)+c\frac{P(X^1_t=1,X^2_t=0)}{P(X^2_t=0)}   & b+c -c\frac{P(X^1_t=1,X^2_t=0)}{P(X^2_t=0)}  \\[2mm]
   0 & 0
\end{array}\right).
\end{eqnarray}
It is easily checked that $A^1(t)=\Lambda ^1(t)$ and $A^2(t)=\Lambda ^2(t)$, as claimed. Consequently, the matrix $\Lambda $ in \eqref{matrgenL} is a weak Markov copula for  $Y^1$ and $Y^2$, but it  is not a strong Markov copula for  $Y^1$ and $Y^2$.

Finally, note that the transition intensities and transition probabilities for $X^1$ do depend on the state of $X^2$:
\begin{itemize}
\item When $X$ is in state $(0,0)$ at some point in time, then, the intensity of transition of $X^1$ from $0$ to $1$ is equal to $a+c$; when  $X$ is in state $(0,1)$ at some point in time, the intensity of transition of $X^1$ from $0$ to $1$ is equal to $a$.
\item When $X$ is in state $(0,0)$ at some point in time, then,  the  transition  probability of $X^1$ from $0$ to $1$ in $t$ units of time is  $$ e^{-bt}(1-e^{-(a+c)t})\frac{a}{a+c}   +1+
  e^{-(a+b+c)t}\left(\frac{a}{a+c}-\frac{c}{b+c}\right)-\frac{a}{a+c}e^{-bt}-\frac{b}{b+c}e^{-at};$$
  when   $X$ is in state $(0,1)$ at some point in time, the  transition  probability of $X^1$ from $0$ to $1$ in $t$ units of time is
  $$ 1-e^{-at}.$$
\end{itemize}
An analogous observation holds for $X^2$, that is, the transition intensities and transition probabilities for $X^2$ do depend on the state of $X^1$.
\eex

\bex\lab{ex3} Here we give an example of a bivariate Markov chain which is not weakly Markovian consistent.

Let us consider the bivariate Markov chain $X=(X^1,X^2)$ on the state space
\begin{equation*}
E=\left\{ {\left( 0,0\right) },{\left(
0,1\right) },{\left( 1,0\right) },{\left( 1,1\right) }\right\}
\end{equation*}
generated by the matrix
\begin{eqnarray}\label{matrgen3}
A=\left(
\begin{array}{cccccccc}
  -(a+b+c)  & b  & a   &
  c  \\[2mm]
  0 & -(d+e)  & d & e   \\[2mm]
  0 & 0 & -f & f  \\[2mm]
  0 & 0 & g & -g
\end{array}\right).
\end{eqnarray}
We denote by  $H^2_{0,1}$ the process that counts the number of jumps of the component $X^2$ from state $0$ to state $1.$ The $\ff^X$-intensity of such jumps is
\be
\lab{inth012} \I_\set{X^1_t=0,X^2_t=0}(b+c)+ \I_\set{X^1_t=1,X^2_t=0}f,
\ee
so the optional projection of this intensity on $\ff^{X^2}$ has the form
\be
\lab{optproj} (b+c)\P(X^1_t=0,X^2_t=0|\F^{X^2}_t)+f\P(X^1_t=1,X^2_t=0|\F^{X^2}_t).
\ee
Since $\{X^2_t=0,X^2_{t/2}=1\} \subseteq \{X^1_{t}=1\}$, on the set $\{X^2_t=0,X^1_{t/2}=1\}$ we have
\be
\P(X^1_t=0,X^2_t=0|X^2_t=0,X^2_{t/2}=1)=0,\quad \P(X^1_t=1,X^2_t=0|X^2_t=0,X^2_{t/2}=1)=1.
\ee
Therefore the above optional projection, on the set $\{X^2_t=0,X^2_{t/2}=1\}$, is equal to
\be
f \P(X^1_t=1,X^2_t=0|X^2_t=0,X^2_{t/2}=1)=f.
\ee
However, on $\{X^2_t=0\}$ the above optional projection is equal to
\beq
&&(b+c)\P(X^1_t=0,X^2_t=0|X^2_t=0)+f\P(X^1_t=1,X^2_t=0|X^2_t=0)\\ &&\quad \quad \quad = (b+c -f ) \P(X^1_t=0,X^2_t=0|X^2_t=0) + f.
\eeq
Assuming that the process $X$ starts from $(0,0)$ at time $t=0$, it can be shown that $\P(X^1_t=0,X^2_t=0|X^2_t=0) > 0$.
 Verification of this is straightforward, but computationally intensive, and can be obtained from the authors on request. Thus, if $b+c \neq f$, then the optional projection on $\F^{X^2}_t$ of the $\ff^X$ intensity of $H^2_{0,1}$ depends on the trajectory of $X^2$ until time $t$, and not just on the state of $X^2$ at time $t$. Thus, $X^2$ is not Markovian in its own filtration. It is obviously not Markovian in the filtration of the entire process $X$ either.
\eex

\end{document}